\documentclass{amsart}
\usepackage{amssymb, amsfonts}
\usepackage[v2]{xy}

\newtheorem{thm}{Theorem}[section]
\newtheorem{cor}[thm]{Corollary}
\newtheorem{prop}[thm]{Proposition}
\newtheorem{lem}[thm]{Lemma}

\theoremstyle{definition}
\newtheorem{defn}[thm]{Definition}

\newtheorem{exmp}[thm]{Example}

\theoremstyle{remark}
\newtheorem{rem}[thm]{Remark}

\DeclareFontFamily{OMS}{rsfs}{\skewchar\font'60}
\DeclareFontShape{OMS}{rsfs}{m}{n}{<-5>rsfs5 <5-7>rsfs7 <7->rsfs10 }{}
\DeclareSymbolFont{rsfs}{OMS}{rsfs}{m}{n}
\DeclareSymbolFontAlphabet{\scr}{rsfs}

\let\overto\xrightarrow

\newcommand{\sA}{\scr{A}}
\newcommand{\sB}{\scr{B}}
\newcommand{\sC}{\scr{C}}
\newcommand{\sD}{\scr{D}}

\newcommand{\sG}{\scr{G}}

\newcommand{\sO}{\scr{O}}

\newcommand{\al}{\alpha}
\newcommand{\be}{\beta}
\newcommand{\ga}{\gamma}
\newcommand{\de}{\delta}
\newcommand{\epz}{\varepsilon}
\newcommand{\ph}{\phi}

\newcommand{\et}{\eta}

\newcommand{\rh}{\rho}
\newcommand{\si}{\sigma}
\newcommand{\ta}{\tau}

\newcommand{\ps}{\psi}
\newcommand{\ze}{\zeta}
\newcommand{\om}{\omega}

\newcommand{\SI}{\Sigma}

\newcommand{\com}{\circ}     
\newcommand{\iso}{\cong}     
\newcommand{\ten}{\otimes}   

\newcommand{\rtarr}{\longrightarrow}

\def\quickop#1{\expandafter\newcommand\csname #1\endcsname{\operatorname{#1}}}
\quickop{Hom} \quickop{End} \quickop{Aut} \quickop{Tel} \quickop{Mic} 
\quickop{Ext} \quickop{Tor} \quickop{Id} \quickop{Coker} \quickop{Ker}
\quickop{Lim} \quickop{Colim} \quickop{Holim} \quickop{Hocolim}
\quickop{id} \quickop{tel} \quickop{mic} \quickop{coker} \quickop{Map}
\quickop{colim} \quickop{holim} \quickop{hocolim} \quickop{im}

\makeatletter
\let\c@equation\c@thm
\makeatother
\numberwithin{equation}{section}

\let\overto\xrightarrow

\thanks{The second and third authors were partially supported by the NSF}

\title{Isomorphisms between left and right adjoints}
\author{H. Fausk, P. Hu, and J.P. May}
\date{\today}

\begin{document}

\begin{abstract}
There are many contexts in algebraic geometry, algebraic topology, and
homological algebra where one encounters a functor that has both a left and 
right adjoint, with the right adjoint being isomorphic to a shift of the left
adjoint specified by an appropriate ``dualizing object''. Typically the left 
adjoint is well understood while the right adjoint is more mysterious, and the 
result identifies the right adjoint in familiar terms. We give a categorical
discussion of such results. One essential point is to differentiate between
the classical framework that arises in algebraic geometry and a deceptively 
similar, but genuinely different, framework that arises in algebraic topology. 
Another is to make clear which parts of the proofs of such results are formal.
The analysis significantly simplifies the proofs of particular cases, as we 
illustrate in a sequel discussing applications to equivariant stable homotopy 
theory.
\end{abstract}

\maketitle

\tableofcontents

We give a categorical discussion of Verdier and Grothendieck isomorphisms
on the one hand and formally analogous results whose proofs involve different 
issues on the other.  Our point is to explain and compare the two contexts and 
to differentiate the formal issues from the substantive issues in each. 
The philosophy goes back to Grothendieck's ``six operations'' formalism.
We fix our general framework, explain what the naive versions of our theorems say, 
and describe which parts of their proofs are 
formal in \S\S1--4. This discussion does not require triangulated categories. 
Its hypotheses and conclusions make sense in general closed symmetric 
monoidal categories, whether or not triangulated. In practice, that means that 
the arguments apply equally well before or after passage to derived categories.

After giving some preliminary results about triangulated 
categories in \S5, we explain the formal theorems comparing left and right
adjoints in \S6. Our ``formal Grothendieck isomorphism theorem'' is an abstraction 
of results of Amnon Neeman, and our ``formal Wirthm\"uller isomorphism theorem'' borrows 
from his ideas.  His paper \cite{Nee} has been influential, and he must be thanked
for catching a mistake in a preliminary version by the third author. We thank Gaunce 
Lewis for discussions of the topological context, and we thank Sasha Beilinson and Madhav 
Nori for making clear that, contrary to our original expectations, the context encountered 
in algebraic topology is not part of the classical context familiar to algebraic geometers.
We also thank Johann Sigurdsson for corrections and emendations. 

\section{The starting point: an adjoint pair $(f^*,f_*)$}

We fix closed symmetric monoidal categories $\sC$ and $\sD$
with respective unit objects $S$ and $T$. We write $\ten$ 
and $\Hom$ for the product and {\em internal} hom functor 
in either category, and we write $X$ (sometimes also $W$) and 
$Y$ (sometimes also $Z$) generically for objects of $\sC$ and 
objects of $\sD$, respectively. We write $\sC(W,X)$ and 
$\sD(Y,Z)$ for the categorical hom sets. We let
$DX = \Hom (X,S)$ denote the dual of $X$.
We let $\text{ev}\colon \Hom(X,W)\ten X\rtarr W$ denote the
evaluation map, that is, the counit of the $(\ten,\Hom)$ adjunction
$$\sC(X\ten X',W)\iso \sC(X,\Hom (X',W)).$$

We also fix a strong symmetric monoidal functor $f^*\colon \sD\rtarr \sC$.
This means that we are given isomorphisms 
\begin{equation}\label{symmon}
f^*T \iso S \ \ \ \text{and}\ \ \ f^*(Y\ten Z)\iso f^*Y\ten f^*Z,
\end{equation}
the second natural, that commute with the associativity, commutativity, and unit 
isomorphisms for $\ten$ in $\sC$ and $\sD$. We assume throughout that $f^*$ has a 
right adjoint $f_{*}$, and we write
$$\epz\colon f^*f_*X \rtarr X\2 \ \ \text{and}\ \ \ \et\colon Y\rtarr f_*f^* Y$$
for the counit and unit of the adjunction. This general context is fixed throughout.

The notation $(f^*,f_*)$ meshes with standard notation in algebraic geometry, where 
one starts with a map $f\colon A\rtarr B$ of spaces or schemes and $f^*$ and $f_*$ are 
pullback and pushforward functors on sheaves. In our generality there need be no 
underlying map ``$f$'' in sight. Some simple illustrative examples are given in \S3.

The assumption that $f^*$ is strong symmetric monoidal has several basic implications.
To begin with, the adjoints of the isomorphism $f^*T\iso S$ and the map
$$\xymatrix@1{
f^*(f_*W\ten f_*X)\iso f^*f_*W\ten f^*f_*X \ar[r]^-{\epz\ten\epz} & W\ten X}\\$$
are maps
\begin{equation}\label{symmon2}
T\rtarr f_*S \ \ \ \text{and}\ \ \ f_*W\ten f_*X\rtarr f_*(W\ten X).
\end{equation}
These are not usually isomorphisms. This means that $f_*$ is {\em lax}\, symmetric monoidal. 

The adjoint of the map
$$\xymatrix@1{
f^*\Hom(Y,Z)\ten f^*Y \iso f^*(\Hom(Y,Z)\ten Y) \ar[r]^-{f^*(\text{ev})}
& f^*Z\\}$$
is a natural map
\begin{equation}\label{al}
\al\colon f^*\Hom(Y,Z)\rtarr \Hom(f^*Y,f^*Z).
\end{equation}
It may or may not be an isomorphism in general, and we say that $f^*$ is 
{\em closed symmetric monoidal}\, if it is. However, the adjoint of the composite map
$$\xymatrix@1{
f^*\Hom(Y,f_*X)\ar[r]^-{\al} &\Hom(f^*Y,f^*f_*X)\ar[rr]^-{\Hom(\id,\epz)} & & \Hom(f^*Y,X)\\}$$
is a natural isomorphism
\begin{equation}\label{recipy2}
\Hom(Y,f_*X) \iso f_*\Hom(f^*Y,X).
\end{equation}
In particular, $\Hom(Y,f_*S)\iso f_*Df^*Y$.
Indeed, we have the following two chains of isomorphisms of functors.
$$ \sD(Z,\Hom(Y,f_*X)) \iso \sD(Z\ten Y,f_*X) \iso  \sC(f^*(Z\ten Y),X)$$
$$ \sD(Z,f_*\Hom(f^*Y,X))  \iso \sC(f^*Z,\Hom(f^*Y,X))\iso \sC(f^*Z\ten f^*Y, X)$$
By the Yoneda lemma and a check of maps, these show immediately that the assumed
isomorphism of functors in (\ref{symmon}) is {\em equivalent\,} to the claimed isomorphism
of functors (\ref{recipy2}). That is, the isomorphism of left adjoints in (\ref{symmon}) 
is adjoint to the isomorphism of right adjoints in (\ref{recipy2}). Systematic 
recognition of such ``conjugate'' pairs of isomorphisms can substitute for quite a 
bit of excess verbiage in the literature. We call this a ``comparison of adjoints'' 
and henceforward leave the details of such arguments to the reader.

Using the isomorphism (\ref{recipy2}), we obtain the following map $\be$, which is 
analogous to both $\al$ and the map of (\ref{symmon2}). Like the latter, it is not
usually an isomorphism.
\begin{equation}\label{be}
\xymatrix@1{
\be\colon f_*\Hom(X,W) \ar[rr]^-{f_*\Hom(\epz,\id)} & & f_*\Hom(f^*f_*X,W)\ar[r]^-{\iso} &\Hom(f_*X,f_*W).\\}
\end{equation}

Using (\ref{symmon2}), we also obtain a natural composite
\begin{equation}\label{recipy1}
\xymatrix@1{
\pi\colon Y\ten f_*X \ar[r]^-{\et\ten\id} & f_*f^*Y\ten f_*X \ar[r] & f_*(f^*Y\ten X).\\}
\end{equation}
Like $\al$, it may or may not be an isomorphism in general. When it is, we say that
the {\em projection formula}\, holds. 

As noted by Lipman \cite[p.\,119]{Lip}, there is already a non-trivial ``coherence problem''
in this general context, the question of determining which compatibility diagrams 
relating the given data necessarily commute. 
An early reference for coherence in closed symmetric monoidal categories is \cite{EK}, and 
the volume \cite{coh} contains several papers on the subject and many references. In particular, 
a paper of G. Lewis in \cite{coh} gives a partial coherence theorem for closed monoidal functors. 
The categorical literature of coherence is relevant to the study of ``compatibilities'' that 
focuses on base change maps and plays an important role in the literature in algebraic geometry 
(e.g. \cite{Lip, Conrad, De2, DeGro, Hart}). A study of that is beyond the scope of this note.
A full categorical coherence theorem is not known and would be highly desirable. 

We illustrate by recording a particular commutative coherence diagram, namely 
\begin{equation}\label{silly}\xymatrix{
f^*DY\ten f^*Y \ar[r]^-{\iso} \ar[d]_{\al\ten\id} & f^*(DY\ten Y) \ar[r]^-{f^*(ev)} 
& f^*T \ar[d]^{\iso}\\
Df^*Y\ten f^*Y \ar[rr]_-{ev} & & S. \\}
\end{equation}
We shall need a consequence of this diagram. There is a natural map
$$\nu\colon DX\ten W\rtarr \Hom(X,W),$$
namely the adjoint of 
$$DX\ten W\ten X \iso DX\ten X\ten W \overto{ev\ten\id} S\ten W\iso W.$$
The commutativity of the diagram (\ref{silly}) implies the commutativity of the diagram
\begin{equation}\label{sillier}\xymatrix{
f^*DY\ten f^*Z \ar[r]^-{\iso}\ar[d]_{\al\ten\id} &  f^*(DY\ten Z) \ar[r]^-{f^*\nu} 
& f^*\Hom(Y,Z) \ar[d]^{\al}\\
Df^*Y\ten f^*Z \ar[rr]_{\nu} & & \Hom(f^*Y,f^*Z).\\}
\end{equation}

We assume familiarity with the theory of ``dualizable'' (alias ``strongly dualizable'' 
or ``finite'') objects; see \cite{May1} for a recent exposition. The defining property
is that $\nu\colon DX\ten X \rtarr \Hom(X,X)$ is an isomorphism. It follows that $\nu$ 
is an isomorphism if either $X$ or $W$ is dualizable. It also follows that the natural 
map $\rh: X\rtarr DDX$ is an isomorphism, but the converse fails in general. When $X'$ is 
dualizable, we have the duality adjunction
\begin{equation}\label{dualad}
\sC (X\ten X',X'') \iso \sC (X, DX'\ten X'').
\end{equation}

As observed in \cite[III.1.9]{LMS}, (\ref{symmon}) and the definitions 
imply the following result.

\begin{prop}\label{easyLMS}
If $Y\in \sD$ is dualizable, then $DY$, $f^*Y$, and $Df^*Y$ are dualizable and the map 
$\al$ of (\ref{al}) restricts to an isomorphism 
\begin{equation}\label{D1}
f^*DY \iso Df^*Y.
\end{equation}
\end{prop}

This implies that $\al$ and $\pi$ are often isomorphisms for formal reasons.

\begin{prop}\label{alpi} If $Y\in\sD$ is dualizable, then
$$\al\colon f^*\Hom(Y,Z)\rtarr \Hom(f^*Y,f^*Z)\ \ \ \text{and} \ \ \ 
\pi\colon Y\ten f_*X \rtarr f_*(f^*Y\ten X)$$
are isomorphisms for all objects $X\in \sC$ and $Z\in\sD$.
Thus, if all objects of $\sD$ are dualizable, then $f^*$ is closed symmetric monoidal and the 
projection formula holds. 
\end{prop}
\begin{proof} For the first statement, $\al$ coincides with the composite 
$$f^*\Hom(Y,Z) \iso f^*(DY\ten Z)\iso f^*DY\ten f^*Z\iso Df^*Y\ten f^*Z\iso \Hom(f^*Y,f^*Z).$$
For the second statement, $\pi$ induces the isomorphism of represented functors
$$\sD(Z,Y\ten f_*X) \iso \sD(Z\ten DY,f_*X) \iso  \sC(f^*(Z\ten DY), X) \iso \sC(f^*Z\ten f^*DY, X)$$
\vspace{-4mm}
$$\hspace*{12mm} \iso  \sC(f^*Z\ten Df^*Y,X)  
\iso  \sC(f^*Z,f^*Y\ten X)\iso \sD(Z,f_*(f^*Y\ten X)). \qed$$
\renewcommand{\qed}{}
\end{proof}

\section{The general context: adjoint pairs $(f^*,f_*)$ and $(f_!,f^!)$}

In addition to the adjoint pair $(f^*,f_*)$ of the previous section, we here
assume given a second adjoint pair $(f_!,f^!)$ relating $\sC$ and $\sD$, with
$f_!\colon \sC\rtarr \sD$ being the left adjoint. We write 
$$\si\colon f_{!}f^{!} Y\rtarr Y \ \ \text{and}\ \ \ze\colon X\rtarr f^{!}f_{!} X$$
for the counit and unit of the second adjunction. 

The adjunction $\sD(Y,f_*X)\iso \sC(f^*Y,X)$ can be recovered from the more general
``internal $\Hom$ adjunction'' $\Hom(Y,f_*X) \iso f_*\Hom(f^*Y,X)$
of (\ref{recipy2}) by applying the functor $\sD(T,-)$ and using the assumption that 
$f^*T\iso S$. Analogously, it is natural to hope that the adjunction 
$\sD(f_!X,Y)\iso\sC(X,f^!Y)$ can be recovered by applying the functor $\sD(T,-)$ 
to a similar internal $\Hom$ adjunction 
$$\Hom(f_!X,Y)\iso f_*\Hom(X,f^!Y).$$

However, unlike (\ref{recipy2}), such an adjunction does not follow formally from
our hypotheses. Motivated by different specializations of the general context, we
consider two triads of basic natural maps that we might ask for relating our four 
functors. For the first triad, we might ask for either of the following two duality
maps, the first of which is a comparison map for the desired internal $\Hom$ adjunction. 
\begin{equation}\label{Ver}
\ga\colon f_*\Hom(X,f^!Y) \rtarr \Hom(f_!X,Y).
\end{equation}
\begin{equation}\label{Ver'}
{\de}\colon \Hom(f^*Y,f^!Z) \rtarr f^!\Hom(Y,Z).
\end{equation}
We might also ask for a {\em projection formula map}
\begin{equation}\label{proj2}
\hat{\pi}\colon Y\ten f_!X \rtarr f_!(f^*Y\ten X),
\end{equation}
which should be thought of as a generalized analogue of the map $\pi$ of (\ref{recipy1}).
These three maps are not formal consequences of the given adjunctions, but rather 
must be constructed by hand. However, it suffices to construct any one of them. 

\begin{prop}\label{three} 
Suppose given any one of the natural maps $\ga$, $\de$, and $\hat{\pi}$. Then it 
determines the other two by conjugation. The map $\de$ is an isomorphism for all 
dualizable $Y$ if and only if its conjugate $\hat{\pi}$ is an isomorphism for all 
dualizable $Y$. If any one of the three conjugately related maps is a natural 
isomorphism, then so are the other two.
\end{prop}

The second triad results from the first simply by changing the direction of the arrows.
That is, we can ask for natural maps in the following directions.
\begin{equation}\label{WirVer}
\bar{\ga} \colon \Hom(f_!X,Y)\rtarr f_*\Hom(X,f^!Y).
\end{equation}
\begin{equation}\label{WirVer'}
\bar{\de}\colon f^!\Hom(Y,Z)\rtarr \Hom(f^*Y,f^!Z).
\end{equation}
\begin{equation}\label{proj3}
\bar{\pi}\colon f_!(f^*Y\ten X)\rtarr Y\ten f_!X.
\end{equation}
Here $\bar{\de}$ is to be viewed as a generalized analogue of the map $\al$ of (\ref{al}).

\begin{prop}\label{three3} 
Suppose given any one of the natural maps $\bar{\ga}$, $\bar{\de}$, and $\bar{\pi}$. Then it 
determines the other two by conjugation. The map $\bar{\de}$ is an isomorphism for all dualizable
$Y$ if and only if its conjugate $\bar{\pi}$ is an isomorphism for all dualizable $Y$. If any one 
of the three conjugately related maps is a natural isomorphism, then so are 
the other two.
\end{prop}

Of course, when the three maps are isomorphisms, the triads are inverse to each other
and there is no real difference. However, there are two very different interesting 
specializations: we might have $f_! = f_*$, or we might have $f^!=f^*$. The first occurs 
frequently in algebraic geometry, and is familiar. The second occurs in algebraic topology 
and elsewhere, but seems less familiar. With the first specialization, the first triad
of maps arises formally since we can take $\hat{\pi}$ to be the map $\pi$ of (\ref{recipy1}). 
With the second specialization, the second triad arises formally since we can take $\bar{\de}$
to be the map $\al$ of (\ref{al}). Recall the isomorphism (\ref{recipy2}), the map $\be$ of 
(\ref{be}), and Proposition \ref{alpi}.

\begin{prop}\label{threea} 
Suppose $f_!=f_*$. Taking $\hat{\pi}$ to be the projection map $\pi$ of (\ref{recipy1}), 
the conjugate map $\ga$ is the composite
$$\xymatrix@1{
f_*\Hom(X,f^!Y) \ar[r]^-{\be} & \Hom(f_*X,f_*f^!Y) \ar[rr]^-{\Hom(\id,\si)} & &  \Hom(f_*X,Y)\\}$$
and the conjugate map $\de$ is the adjoint of the map
$$\xymatrix@1{
f_*\Hom(f^*Y,f^!Z) \iso \Hom(Y,f_*f^!Z) \ar[rr]^-{\Hom(\id,\si)} & & \Hom(Y,Z).\\}$$
Moreover, $\pi$ and $\de$ are isomorphisms if $Y$ is dualizable.
\end{prop}

When $f^!=f^*$, passage to adjoints from $S\iso f^*T$ and the natural map
$$\xymatrix@1{
W\ten X\ar[r]^-{\ze\ten\ze} & f^*f_!W\ten f^*f_!X \iso f^*(f_!W\ten f_!X)\\}$$
gives maps, not usually isomorphisms, 
\begin{equation}\label{symmon3}
f_!S\rtarr T \ \ \ \text{and}\ \ \ f_!(W\ten X)\rtarr f_!W\ten f_!X.
\end{equation}
This means that $f_!$ is an {\em op-lax}\, symmetric monoidal functor.

\begin{prop}\label{threeb}
Suppose $f^!=f^*$. Taking $\bar{\de}$ to be the map $\al$ of (\ref{al}), the 
conjugate map $\bar{\pi}$ is the composite
$$\xymatrix@1{
f_!(f^{*}Y\ten X) \ar[r] & f_!f^{*}Y\ten f_!X
\ar[r]^-{\si\ten\id} & Y\ten f_!X\\}$$
and the conjugate map $\bar{\ga}$ is the adjoint of the map
$$\xymatrix@1{
f^*\Hom(f_!X,Y)\ar[r]^-{\al} & \Hom(f^*f_!X,f^*Y)\ar[rr]^-{\Hom(\ze,\id)} 
& & \Hom(X,f^*Y).\\}$$
Moreover $\al$ and $\bar{\pi}$ are isomorphisms if $Y$ is dualizable.
\end{prop}

\begin{defn} We introduce names for the different contexts in sight.\\
(i) {\em The Verdier-Grothendieck context:} There is a natural isomorphism $\hat{\pi}$ as in
(\ref{proj2}) (projection formula); taking $\bar{\pi} = \hat{\pi}^{-1}$, there are conjugately
determined natural isomorphisms $\ga = \bar{\ga}^{-1}$, and $\de = \bar{\de}^{-1}$.\\
(ii) {\em The Grothendieck context:} $f_!= f_*$ and the projection formula holds.\\
(iii) {\em The Wirthm\"uller context:} $f^!=f^*$ and $f^*$ is closed symmetric monoidal.
\end{defn}

Thus, in the Grothendieck context, the strong symmetric monoidal functor $f^*$ is 
the left adjoint of a left adjoint. In the Wirthm\"uller context, it is a left and 
a right adjoint.

The Verdier-Grothendieck context encapsulates the properties that hold for 
suitable derived categories $\sC$ and $\sD$ of sheaves over locally compact 
spaces $A$ and $B$ and maps $f:A\rtarr B$; see \cite{Borel, Iv, V}. Here $f_!$ 
is given by pushforward with compact supports. The same abstract context applies 
to suitable derived categories $\sC$ and $\sD$ of complexes of $\sO_A$-modules 
and of $\sO_B$-modules for schemes $A$ and $B$ and maps $f:A\rtarr B$. In either 
context, we have $f_!=f_*$ when the map $f$ is proper. For the scheme theoretic
context, see \cite{De, DeGro, Hart} and, for more recent reworkings and generalizations, 
\cite{Lip, Conrad, Lip2, Nee}. There is a highly non-trivial categorical, more 
precisely $2$-categorical, coherence problem concerning composites of base change 
functors in the Verdier-Grothendieck context. A start on this has been made by 
Voevodsky \cite{De2}. Since his discussion focuses on base change relating quadruples 
$(f^*,f_*,f_!,f^!)$, ignoring $\ten$ and $\Hom$, it is essentially 
disjoint from our discussion. The relevant coherence problem simplifies greatly in 
either the Grothendieck or the Wirthm\"uller context, due to the canonicity of the 
maps in Propositions \ref{threea} and \ref{threeb}. 

We repeat that our categorical results deduce formal conclusions from formal 
hypotheses and therefore work equally well before or after passage to derived 
categories. Much of the work in passing from categories of sheaves to derived 
categories can be viewed as the verification that formal properties in the category 
of sheaves carry over to the same formal properties in derived categories, although
other properties only hold after passage to derived categories.  A paper by Lipman
in \cite{Lip} takes a similarly categorical point of view. 

While the proofs of Propositions \ref{three} and \ref{threea} are formal, 
in the applications to algebraic geometry they require use of {\em unbounded}\, 
derived categories, since otherwise we would not have closed symmetric monoidal categories 
to begin with. These were not available until Spaltenstein's paper \cite{Spalt}, and he 
noticed one of our formal implications \cite[\S6]{Spalt}. Unfortunately, as he makes clear, 
in the classical sheaf context his methods fail to give the $(f_!,f^!)$ adjunction 
for all maps $f$ between locally compact spaces. It seems possible that a model theoretic 
approach to unbounded derived categories would allow one to resolve this problem. In any case, 
a complete reworking of the theory in model theoretical terms would be of considerable value.

In the algebraic geometry setting, smooth maps lead to a context close to the Wirthm\"uller
context, but that is not our motivation. In that context, we think of $f^*$ as a forgetful 
functor which does not alter underlying structure, $f_!$ as a kind of extension of scalars 
functor, and $f_*$ as a kind of ``coextension of scalars'' functor. 

For example, let $f: H\rtarr G$ be an inclusion of a subgroup in a group $G$ and let 
$\sC$ and $\sD$ be the categories of $H$-objects and $G$-objects in some Cartesian 
closed category, such as topological spaces. Let $f^*:\sD\rtarr \sC$ be the evident 
forgetful functor. Certainly  
$$f^*(Y\times Z)\iso f^*Y\times f^*Z.$$
The left and right adjoints of $f^*$ send an $H$-object $X$ to $G\times_H X$ and to
$\text{Map}_H(G,X)$. Clearly $G\times_H(X\times X')$ is not isomorphic to 
$(G\times_H X)\times (G\times_H X')$. Our motivating example is a spectrum level 
analogue of this for which there is a Wirthm\"uller isomorphism theorem \cite{LMS, W}. 
Our formal Wirthm\"uller isomorphism theorem below substantially simplifies its 
proof \cite{M}. One can hope for such a result in any context where group actions and 
triangulated categories mix.

For another example, let $f:A\rtarr B$ be an inclusion of cocommutative Hopf algebras over a 
field $k$ and let $\sC$ and $\sD$ be the categories of $A$-modules and of $B$-modules. These 
are closed symmetric monoidal categories under the functors $\otimes_k$ and $\Hom_k$. Indeed, 
using the coproduct on $A$, we see that if $M$ and $N$ are $A$-modules, then so are $M\ten_kN$ 
and $\Hom_k(M,N)$. The commutativity of $\ten_k$ requires the cocommutativity of $A$. The unit 
object in both $\sC$ and $\sD$ is $k$. Again, if $f^*:\sD\rtarr \sC$ is the evident forgetful 
functor, then
$$f^*(Y\ten_k Z)=f^*Y\ten_k f^*Z.$$
The left and right adjoints of $f^*$ send an $A$-module $X$ to $B\ten_A X$ and to
$\Hom_A(B,X)$, and again $B\ten_A(X\ten_k X')$ is not isomorphic to 
$(B\ten_A X)\otimes_k (B\ten_A X')$.
This example deserves investigation on the level of derived categories.

\section{Isomorphisms in the Verdier--Grothendieck context}

We place ourselves in the Verdier--Grothendieck context in this section.

\begin{defn} For an object $W\in \sC$, define $D_W X = \Hom(X,W)$, the 
{\em $W$-twisted dual of $X$}. Of course, if $X$ or $W$ is dualizable, then
$D_W X\iso DX\ten W$. Let $\rh_W: X\rtarr D_WD_WX$ be the adjoint of the 
evaluation map $D_WX\ten X \rtarr W$. We say that $X$ is 
{\em $W$-reflexive}\/ if $\rh_W$ is an isomorphism.
\end{defn}

Replacing $Y$ by $Z$ in (\ref{Ver}) and letting $W=f^!Z$, the isomorphisms 
$\ga$ and $\de$ take the following form:
\begin{equation}\label{dualform}
f_*D_WX \iso D_Z f_!X\ \ \text{and} \ \ D_Wf^*Y \iso f^!D_ZY.
\end{equation}
This change of notation and comparison with the classical context of algebraic
geometry explains why we think of $\ga$ and $\de$ as duality maps. If $f_!X$ is
$Z$-reflexive, the first isomorphism implies that 
\begin{equation}\label{onesort}
f_!X \iso D_Z f_*D_W X.
\end{equation}
If $Y$ is isomorphic to $D_ZY'$ for some $Z$-reflexive object $Y'$, the second
isomorphism implies that
\begin{equation}\label{twosort}
f^!Y \iso D_W f^*D_Z Y.
\end{equation}
These observations and the classical context suggest the following definition.

\begin{defn}\label{dobject} 
A {\em dualizing object} for a full subcategory $\sC_0$ of
$\sC$ is an object $W$ of $\sC$ such that if $X\in \sC_0$, then $D_WX$ is 
in $\sC_0$ and $X$ is $W$-reflexive. Thus $D_W$ specifies an auto--duality
of the category $\sC_0$.
\end{defn}

\begin{rem}\label{dcontext}
In algebraic geometry, we often encounter canonical subcategories 
$\sC_0\subset \sC$ and $\sD_0\subset \sD$ such that $f_!\sC_0\subset \sD_0$
and $f^!\sD_0\subset \sC_0$ together with a dualizing object $Z$ for $\sD_0$ 
such that $W=f^!Z$ is a dualizing object for $\sC_0$. In such contexts, 
(\ref{onesort}) and (\ref{twosort}) express $f_!$ on $\sC_0$ and $f^!$ on 
$\sD_0$ in terms of $f_*$ and $f^*$.
\end{rem}

For any objects $Y$ and $Z$ of $\sD$, the adjoint of the map
$$\xymatrix@1{
f_!(f^*Y\ten f^{!}Z)\iso Y\ten f_!f^{!}Z \ar[r]^-{\id\ten \si}&  Y\ten Z\\}$$
is a natural map 
\begin{equation}\label{YYYV}
\ph\colon f^*Y\ten f^{!}Z \rtarr f^{!}(Y\ten Z).
\end{equation}
It specializes to
\begin{equation}\label{phiV}
\ph\colon f^*Y\ten f^{!}T \rtarr f^{!}Y,
\end{equation}
which of course compares a right adjoint to a shift of a left adjoint. 
A {\em Verdier--Grothendieck isomorphism theorem} asserts that the map 
$\ph$ is an isomorphism; in the context of sheaves over spaces, such a 
result was announced by Verdier in \cite[\S5]{V}. The following observation,
abstracts a result of Neeman \cite[5.4]{Nee}. In it, we only assume the 
projection formula for dualizable $Y$. 

\begin{prop}\label{dualGr} 
The map $\ph\colon f^*Y\ten f^{!}Z \rtarr f^{!}(Y\ten Z)$ is 
an isomorphism for all objects $Z$ and all dualizable objects $Y$.
\end{prop}
\begin{proof}
Using Proposition \ref{easyLMS}, the projection formula, duality 
adjunctions (\ref{dualad}), and the $(f_!,f^{!})$ adjunction, 
we obtain isomorphisms
$$ \sC(X,f^*Y\ten f^{!}Z) \iso \sC(f^*DY\ten X,f^{!}Z) \iso  \sD(f_!(f^*DY\ten X),Z)$$
\vspace{-4mm}
$$\hspace*{12mm}   \iso \sD(DY\ten f_!X, Z)  \iso  \sD(f_!X, Y\ten Z)  
\iso  \sC(X,f^{!}(Y\ten Z)).$$
Diagram chasing shows that the composite isomorphism is induced by $\ph$.
\end{proof}

It is natural to ask when $\ph$ is an isomorphism in general, and we shall return to 
that question in the context of triangulated categories. Of course, this discussion 
specializes and remains interesting in the Grothendieck context $f_!=f_*$. 

We give some elementary examples of the Verdier--Grothendieck context.

\begin{exmp}\label{trivex}
An example of the Verdier--Grothendieck context is already available 
with $\sC=\sD$ and $f^* = f_* = \Id$. Fix an object $C$ of $\sC$ and set
$$ f_!X = X\ten C \ \ \text{and} \ \ f^!(Y) = \Hom(C,Y).$$
The projection formula $f_!(f^*Y\ten Z)\iso Y\ten f_!Z$ is 
the associativity isomorphism
$$(Y\ten Z)\ten C\iso Y\ten (Z\ten C).$$
The map $\ph\colon f^*Y\ten f^!Z\rtarr f^!(Y\ten Z)$ is the canonical map
$$\nu\colon Y\ten \Hom(C,Z) \rtarr \Hom(C,Y\ten Z).$$ 
It is an isomorphism if $Y$ is dualizable, and it is an isomorphism for all
$Y$ if and only if $C$ is dualizable.
\end{exmp} 

The shift of an adjunction by an object of $\sC$ used in the previous example 
generalizes to give a shift of any Verdier-Grothendieck context by an object 
of $\sC$. 

\begin{defn}
For an adjoint pair $(f_!,f^!)$ and an object $C\in \sC$, define the twisted 
adjoint pair $(f^C_!,f^!_C)$ by
\begin{equation}\label{shift1}
f_{!}^C(X) = f_{!}(X\ten C) \ \ \text{and}\ \ f^{!}_CY = \Hom(C,f^{!}Y).
\end{equation}
\end{defn}

\begin{prop}
If $(f^*,f_*)$ and $(f_!,f!)$ are in the Verdier-Grothendieck context, then so
are $(f^*,f_*)$ and $(f_!^C,f!_C)$. 
\end{prop}
\begin{proof} The isomorphism $\hat{\pi}$ of (\ref{proj2}) shifts to a corresponding 
isomorphism $\hat{\pi}_C$.
\end{proof}

We also give a simple example of the context of Definition \ref{dobject}.
Recall that dualizable objects are $S$-reflexive, but not conversely in general. 
The following observation parallels part of a standard characterization of 
``dualizing complexes'' \cite[V.2.1]{Hart}. Let $d\sC$ denote the full subcategory of 
dualizable objects of $\sC$.

\begin{prop} $S$ is $W$-reflexive if and only if all $X\in d\sC$ are $W$-reflexive.
\end{prop}
\begin{proof} Since $S$ is dualizable, the backwards implication is trivial. 
Assume that $S$ is $W$-reflexive. Since $W\iso D_WS$, $\Hom(W,W)=D_WW\iso D_WD_W S$. 
In any closed symmetric monoidal category, such as $\sC$, we have a natural isomorphism
$$\Hom(X\ten X',X'')\iso \Hom(X,\Hom(X',X'')),$$
where $X$, $X'$, and $X''$ are arbitrary objects. When $X$ is dualizable, 
$$\nu\colon DX\ten X'\rtarr \Hom(X,X')$$ 
is an isomorphism for any object $X'$. Therefore
$$D_WD_WX \iso \Hom(DX\ten W,W)\iso \Hom(DX,\Hom(W,W))\iso DDX\ten D_WD_WS.$$
Identifying $X$ with $X\ten S$, is easy to check that $\rh_W$ corresponds under 
this isomorphism to $\rh_S\ten \rh_W$. The conclusion follows.
\end{proof}

\begin{cor} Let $W$ be dualizable. Then the following are equivalent.
\begin{enumerate}
\item[(i)] $W$ is a dualizing object for $d\sC$.
\item[(ii)] $S$ is $W$-reflexive.
\item[(iii)] $W$ is invertible.
\item[(iv)] $D_W\colon d\sC^{op}\rtarr d\sC$ is an auto--duality of $d\sC$.
\end{enumerate}
\end{cor}
\begin{proof} If $X$ is dualizable, then $D_WX\iso DX\ten W$ is dualizable.
The proposition shows that (i) and (ii) are equivalent, and it is clear
that (iii) and (iv) are equivalent. Since $W$ is dualizable, 
$D_WD_WS\iso \Hom(W,W)\iso W\ten DW$, 
with $\rh_W$ corresponding to the coevaluation map $coev: S\rtarr W\ten DW$. 
By \cite[2.9]{May1}, $W$ is invertible if and only if $coev$ is an isomorphism. 
Therefore (ii) and (iii) are equivalent.
\end{proof}

Finally, we have a shift comparison of Grothendieck and Wirthm\"uller contexts.

\begin{rem} Start in the Grothendieck context, so that $f_!=f_*$,  and assume that 
the map $\phi\colon f^*Y\ten f^!T\rtarr f^!Y$ of (\ref{phiV}) is an isomorphism.
Assume further that $f^!T$ is invertible and let $C=Df^!T$. Define a new functor
$f_!$ by $f_!X = f_*(X\ten DC)$. Then $f_!$ is {\em left}\/ adjoint to $f^*$. 
Replacing $X$ by $X\ten C$, we see that 
$$f_*X \iso f_!(X\ten C).$$
In the next section, we shall consider isomorphisms of this general form in the 
Wirthm\"uller context.  Conversely, start in the Wirthm\"uller context, so that 
$f^!=f^*$, and assume given a $C$ such that $f_*S\iso f_!C$ and the map 
$\om\colon f_*X\rtarr f_!(X\ten C)$ of (\ref{omega}) below is an isomorphism. 
Define a new functor $f^!$ by $f^!Y=\Hom(C,f^*Y)$ and note that $f^!T \iso DC$. 
Then $f^!$ is {\em right}\/ adjoint to $f_*$. If either $C$ or $Y$ is dualizable, 
then $\Hom(C,f^*Y)\iso f^*Y\ten DC$ and thus $f^*Y\ten f^!T \iso f^!Y$, which is 
an isomorphism of the same form as in the Grothendieck context.
\end{rem}

\section{The Wirthm\"uller isomorphism}

We place ourselves in the Wirthm\"uller context in this section, with 
$f^!=f^*$. Here the specialization of the Verdier--Grothendieck isomorphism is of no interest.
In fact, $\ph$ reduces to the originally assumed isomorphism (\ref{symmon}). However, 
there is now a candidate for an isomorphism between the right adjoint $f_*$ of $f^*$ 
and a shift of the left adjoint $f_!$. This is not motivated by duality questions, and 
it can already fail on dualizable objects. We assume in addition to the isomorphisms 
$\al = \bar{\de}$, hence $\bar{\pi}$ and $\bar{\ga}$, that we are given an object $C\in\sC$ 
together with an isomorphism
\begin{equation}\label{Dobject}
f_*S\iso f_!C. 
\end{equation}
Observe that the isomorphism $\bar{\ga}$ specializes to an isomorphism
\begin{equation}\label{D3} 
Df_!X \iso f_*DX.
\end{equation}
Taking $X=S$ in (\ref{D3}) and using that $DS\iso S$, we see that (\ref{Dobject}) is 
equivalent to 
\begin{equation}\label{Dobjalt}
Df_{!} S\iso f_{!}C.
\end{equation}
This version is the one most naturally encountered in applications, since it 
makes no reference to the right adjoint $f_*$ that we seek to understand.
In practice, $f_!S$ is dualizable and $C$ is dualizable or even invertible.
It is a curious feature of our discussion that it does not require such hypotheses.

Replacing $C$ by $S\ten C$ in (\ref{Dobject}), it is reasonable to hope that it 
continues to hold with $S$ replaced by a general $X$. That is, we can hope for a 
natural isomorphism 
\begin{equation}\label{Wirth}
f_*X\iso f_{\sharp}X, \ \ \text{where} \ \ f_{\sharp}X \equiv f_!(X\ten C).
\end{equation}
Note that we twist by $C$ before applying $f_!$. We shall shortly define a particular
natural map $\om\colon f_*X\rtarr f_{\sharp}X$. A {\em Wirthm\"uller isomorphism theorem}\/
asserts that $\om$ is an isomorphism. We shall show that if $f_!S$ is dualizable and $X$ is a 
retract of some $f^*Y$, then $\om$ is an isomorphism. However, even for dualizable $X$, 
$\om$ need {\em not}\, be an isomorphism in general. 
A counterexample is given in the sequel \cite{M}.  We shall also give a categorical criterion 
for $\om$ to be an isomorphism for a particular object $X$. An application is also given in \cite{M}. 

Using the map $T\rtarr f_*S$ of (\ref{symmon2}), the assumed isomorphism 
$f_*S\iso f_{!}C$ gives rise to maps
$$\xymatrix@1{
\ta\colon T\rtarr f_*S \iso f_!C\\}
$$
and
$$\xymatrix@1{
\xi\colon f^*f_!C \iso f^*f_*S \ar[r]^-{\epz} & S\\}
$$
such that
$$\xi\com f^*\ta = \id\colon S\rtarr S.$$
Using the alternative defining property (\ref{Dobjalt}) of $C$, we can obtain alternative 
descriptions of these maps that avoid reference to the functor $f_*$ we seek to understand.

\begin{lem}\label{altalt} 
The maps $\ta$ and $\xi$ coincide with the maps
$$\xymatrix@1{
T\iso DT \ar[r]^-{D\si} & Df_!f^*T \iso Df_!S \iso f_!C\\}
$$
and
$$\xymatrix@1{
f^*f_!C \iso f^*Df_!S\iso Df^*f_!S
\ar[r]^-{D\ze} & DS \iso S.\\}
$$
\end{lem}
\begin{proof} The isomorphism $Df_!S\iso f_!C$ used in the 
displays is the composite of the given isomorphism (\ref{Dobject}) 
and the special case (\ref{D3}) of the isomorphism $\bar{\ga}$. The proofs 
are diagram chases that use the naturality of $\et$ 
and $\epz$, the triangular identities for the $(f_!,f^*)$ adjunction, and the
description of $\bar{\ga}$ in Proposition \ref{threeb}.
\end{proof}

Using the isomorphism $(\ref{proj3})$, we extend $\ta$ to the natural map 
\begin{equation}\label{tau}
\xymatrix@1{
\ta\colon Y\iso Y\ten T \ar[r]^-{\id\ten\ta} & Y\ten f_!C
\iso f_!(f^*Y\ten C) = f_{\sharp}f^*Y.\\}
\end{equation}
Specializing to $Y=f_*X$, we obtain the desired comparison map $\om$ as the composite
\begin{equation}\label{omega}
\xymatrix@1{
\om\colon f_*X \ar[r]^-{\ta} & f_{\sharp}f^*f_*X 
\ar[r]^-{f_{\sharp}\epz} & f_{\sharp}X.\\}
\end{equation}
An easy diagram chase using the triangular identity $\epz\com f^*\et = \id$ shows that 
\begin{equation}\label{tauom}
\om\com \et = \ta\colon Y \rtarr f_{\sharp}f^* Y.
\end{equation}
If $\om$ is an isomorphism, then $\ta$ must be the unit of the resulting $(f^*,f_{\sharp})$ 
adjunction.

Similarly, using (\ref{symmon}) and (\ref{proj3}), we extend $\xi$ to the natural map 
\begin{equation}\label{xipa}
\xymatrix@1{
\xi\colon f^*f_{\sharp}f^*Y = f^*f_!(f^*Y\ten C) 
\iso f^*Y\ten f^*f_!C \ar[r]^-{\id\ten \xi} & f^*Y\ten S\iso f^*Y.\\}
\end{equation}
We view $\xi$ as a partial counit, defined not for all $X$ but only for $X=f^*Y$.
Since $\xi\com f^*\ta =\id\colon S\rtarr S$, it is immediate that
\begin{equation}\label{oneform}
\xi\com f^*\ta =\id\colon f^*Y\rtarr f^*Y,
\end{equation}
which is one of the triangular identities for the desired 
$(f^*,f_{\sharp})$ adjunction. Define
\begin{equation}\label{psi}
\ps\colon f_{\sharp}f^*Y\rtarr f_*f^* Y
\end{equation}
to be the adjoint of $\xi$. The adjoint of the relation (\ref{oneform}) is the 
analogue of (\ref{tauom}):
\begin{equation}\label{omtau}
\ps\com \ta = \et\colon Y\rtarr f_*f^* Y.
\end{equation} 

\begin{prop} If $Y$ or $f_! S$ is dualizable, then
$\om: f_*f^*Y\rtarr f_{\sharp}f^*Y$ is an isomorphism with inverse 
$\ps$. If $\ps$ is an isomorphism for all $Y$, then $f_!S$ is 
dualizable. If $X$ is a retract of some $f^*Y$, where $Y$ or $f_!S$
is dualizable, then $\om: f_*X\rtarr f_{\sharp}X$ is an isomorphism.
\end{prop}
\begin{proof}
With $X=f^*Y$, the first part of the proof of the following result gives that 
$\ps\com\om = \id$, so that $\om = \ps^{-1}$ when $\ps$ is an isomorphism. We 
claim that $\ps$ coincides with the following composite:
$$ f_{\sharp}f^*Y = f_{!}(f^*Y\ten C)\iso Y\ten D(f_{!}S)
\overto{\nu}\Hom(f_{!} S, Y)\iso f_*\Hom(S,f^*Y) = f^*Y.$$
Here the isomorphisms are given by (\ref{proj3}) and (\ref{Dobjalt})
and by (\ref{WirVer}). Since $\nu$ is an isomorphism if $Y$ or $f_{!}S$
is dualizable, the claim implies the first statement. Note that $\ps = f_*\xi\com\et$
and that the isomorphism $\bar{\ga}$ of (\ref{WirVer}) is $f_*\Hom(\ze,\id)\com f_*\al\com \et$.
Using the naturality of $\et$ and the description of $\xi$ in Lemma \ref{altalt}, an easy,
if lengthy, diagram chase shows that the diagram (\ref{sillier}) gives just what is needed 
to check the claim. The second statement is now clear by the definition of dualizability: it
suffices to consider $Y = f_!S$. The last statement follows from the first since a retract 
of an isomorphism is an isomorphism. 
\end{proof}

We extract a criterion for $\om$ to be an isomorphism for a general object $X$ from 
the usual proof of the uniqueness of adjoint functors \cite[p. 85]{Mac}. 
\begin{prop}\label{onebyone} 
If there is a map $\xi\colon  f^*f_{\sharp} X= f^*f_{!}(X\ten C)\rtarr X$ such that 

\begin{equation}\label{keypt}
f_{\sharp}\xi\com \ta = \id\colon f_{\sharp}X\rtarr f_{\sharp}X
\end{equation}
and the following (partial naturality) diagram commutes, then 
$\om\colon f_*X\rtarr f_{\sharp}X$ is an 
isomorphism with inverse the adjoint $\ps$ of $\xi$.
\begin{equation}\label{natdiag}
\xymatrix{
f^*f_{\sharp}f^*f_* X \ar[r]^-{\xi} \ar[d]_{f^*f_{\sharp}\epz} & f^*f_* X \ar[d]^{\epz} \\
f^*f_{\sharp} X \ar[r]_{\xi} & X \\}
\end{equation}
Moreover, (\ref{keypt}) holds if and only if the following diagram commutes.
\begin{equation}\label{altkey}
\xymatrix{
X\ten C \ar[r]^-{\ze} \ar[d]_{\ze} & f^*f_{!}(X\ten C) \\
f^*f_{!}(X\ten C) \ar[r]_-{f^*\ta} 
& f^*f_{!}(f^*f_{!}(X\ten C)\ten C) 
\ar[u]_{f^*f_{!}(\xi\ten\id)} \\}
\end{equation}
\end{prop}
\begin{proof} In the diagram (\ref{natdiag}), the top map $\xi$ is given by (\ref{xipa}). 
The diagram and the relation $\xi\com f^*\ta=\id$ of (\ref{oneform}) easily imply the relation 
$\xi\com f^*\om = \epz$, which is complementary to the defining relation $\epz\com f^*\ps=\xi$ 
for the adjoint $\ps$. Passage to adjoints gives that $\ps\com \om = \id$. The following
diagram commutes by  (\ref{tauom}), the triangular identity $f_*\epz\com\et =\id$, the naturality 
of $\et$ and $\om$, and the fact that $\ps$ is adjoint to $\xi$.
It gives that $\om\com \ps = f_{\sharp}\xi\com \ta = \id$. 

$$\xymatrix{
f_{\sharp}X \ar[rrr]^-{\ta} \ar[dr]^{\et} \ar[ddd]_{\ps} & & & f_{\sharp}f^*f_{\sharp}X \ar[ddd]^{f_{\sharp}\xi} \ar@{=}[dl] \\
& f_*f^*f_{\sharp}X \ar[r]^-{\om} \ar[d]_{f_*f^*\ps} & f_{\sharp}f^*f_{\sharp}X \ar[d]^{f_{\sharp}f^*\ps} &  \\
& f_*f^*f_*X \ar[r]^-{\om} \ar[d]^{f_*\epz}  & f_{\sharp}f^*f_*X \ar[dr]^{f_{\sharp}\epz}& \\
f_*X \ar@{=}[r] \ar[ur]^{\et} & f_*X \ar[rr]_-{\om} & & f_{\sharp}X. \\}$$
The last statement is clear by adjunction.
\end{proof}

\begin{rem} The map $\om$ can be generalized to the Verdier--Grothendieck context. 
For that, we assume given an object $W$ of $\sC$ such that 
$$f_{!}C\iso Df_{!}f^{!}T;$$
compare (\ref{Dobjalt}). As in Lemma \ref{altalt}, we then have the map 
$$\xymatrix@1{
\ta\colon T\iso DT \ar[r]^-{D\si} & Df_!f^!T \iso f_!C.\\}$$
This allows us to define the comparison map
$$\xymatrix@1{
\om\colon f_*X \iso f_*X\ten T\ar[r]^-{\id\ten\ta} & f_*X\ten f_!W \iso f_!(f^*f_*X\ten C)
\ar[r]^-{f_!(\epz\ten\id)} & f_{!}(X\ten C.\\}$$
\end{rem}

A study of when this map $\om$ is an isomorphism might be of interest, but 
we have no applications in mind. We illustrate the idea in the context of 
Example \ref{trivex}. 

\begin{exmp} Returning to Example \ref{trivex}, we seek an object $C'$ of $\sC$ such that 
$f_{!}C' \iso D(f_{!}f^{!}S)$, which is 
$$C'\ten C\iso D(DC\ten C).$$
If $C$ is dualizable, then the right side is isomorphic to $C\ten DC\iso DC\ten C$ and 
we can take $C'=DC$. Here the map 
$$\om\colon X = f_*X\rtarr f_{!}(X\ten DC) = X\ten DC\ten C$$
turns out to be $\id\ten ({\ga}\com coev)$, where $coev\colon S\rtarr C\ten DC$ 
is the coevaluation map of the duality adjunction (\ref{dualad}) and ${\ga}$ is
the commutativity isomorphism for $\ten$. We conclude (e.g., by \cite[2.9]{May1}) 
that $\om$ is an isomorphism if and only if $C$ is invertible.
\end{exmp}

\section{Preliminaries on triangulated categories} 

We now go beyond the hypotheses of \S\S1--4 to the 
triangulated category situations that arise in practice.
We assume that $\sC$ and $\sD$ are triangulated and that the functors $(-)\ten X$ and $f^*$ are exact 
(or triangulated). This means that they are additive, commute up to isomorphism with $\SI$, and preserve 
distinguished triangles. For $(-)\ten X$, this is a small part of the appropriate compatibility 
conditions that relate distinguished triangles to $\ten$ and $\Hom$ in well-behaved triangulated 
closed symmetric monoidal categories; see \cite{May2} for a discussion of this, as well as for basic 
observations about what triangulated categories really are: the standard axiom system is redundant 
and unnecessarily obscure. We record the following easily proven observation relating adjoints to 
exactness (see for example \cite[3.9]{Nee2}).

\begin{lem}\label{adexact}
Let  $F\colon \sA\rtarr \sB$ and $G:\sB\rtarr \sA$ be left and right adjoint functors  
between triangulated categories. Then $F$ is exact if and only if $G$ is exact.
\end{lem}

We also record the following definitions (see for example \cite{HPS, Nee}).

\begin{defn} 
A full subcategory $\sB$ of a triangulated category $\sC$ is {\em thick}
if any retract of an object of $\sB$ is in $\sB$ and if the third object of a 
distinguished triangle with two objects in $\sB$ is also in $\sB$. The 
category $\sB$ is {\em localizing} if it is thick and closed under coproducts.
The smallest thick (respectively, localizing) subcategory of $\sC$ that contains a 
set of objects $\sG$ is called the thick (respectively, localizing) subcategory 
generated by $\sG$.
\end{defn}

\begin{defn}\label{compact} 
An object $X$ of an additive category $\sA$ is {\em compact}, or {\em small}, 
if the functor $\sA(X,-)$ converts coproducts to direct sums. The category $\sA$ is 
{\em compactly generated} if it has arbitrary coproducts and has a set $\sG$ of compact 
objects that detects isomorphisms, in the sense that a map $f$ in $\sA$ is an isomorphism 
if and only if $\sA(X,f)$ is an isomorphism for all $X\in\sG$. When $\sA$ is symmetric
monoidal, we require its unit object to be compact; thus it can be included in the 
set $\sG$.
\end{defn}

In the triangulated case, this is equivalent to Neeman's definition \cite[1.7]{Nee}.
With our version, we have the following generalization of a result of his \cite[5.1]{Nee}.

\begin{lem}\label{preserve} Let $\sA$ be a compactly generated additive
category with generating set $\sG$ and let $\sB$ be any additive category. 
Let $F\colon \sA\rtarr \sB$ be an additive functor with right adjoint $G$. 
If $G$ preserves coproducts, then $F$ preserves compact objects. Conversely, 
if $F(X)$ is compact for $X\in \sG$, then $G$ preserves coproducts.
\end{lem}
\begin{proof}
Let $X\in\sA$ and let $\{Y_i\}$ be a set of objects of $\sB$. Then the
evident map $f:\amalg G(Y_i)\rtarr G(\amalg Y_i)$ induces a map
$$f_*\colon \sA(X,\amalg G(Y_i)) \rtarr \sA(X,G(\amalg Y_i)).$$
If $X$ is compact and $f_*$ is an isomorphism, then, by adjunction and compactness,
it induces an isomorphism 
$$\amalg \sB(F(X), Y_i) \rtarr \sB(F(X),\amalg Y_i),$$
which shows that $F(X)$ is compact. Conversely, if $X$ and $F(X)$ are both compact,
then $f_*$ corresponds under adjunction to the identity map of $\amalg \sB(F(X),Y_i)$
and is therefore an isomorphism. Restricting to $X\in\sG$, it follows from 
Definition \ref{compact} that $f$ is an isomorphism.
\end{proof}

While this result is elementary, it is fundamental to the applications. We generally
have much better understanding of left adjoints, so that the compactness criterion
is verifiable, but it is the preservation of coproducts by right adjoints that is
required in all of the formal proofs.

Returning to triangulated categories, we justify the term ``generating set'' by 
the following result. Its first part is \cite[3.2]{Nee}, and its second part is 
\cite[2.1.3(d)]{HPS}.

\begin{prop}\label{thickgen}
Let $\sA$ be a compactly generated triangulated category with generating set $\sG$.
Then the localizing subcategory generated by $\sG$ is $\sA$ itself. If the objects
of $\sG$ are dualizable, then the thick subcategory generated by $\sG$ is the full 
subcategory of dualizable objects in $\sA$, and an object is dualizable if and only
if it is compact.
\end{prop}

The following standard observation works in tandem with the previous result.

\begin{prop}\label{wrap}
Let $F, F'\colon \sA\rtarr \sB$ be exact functors between triangulated categories 
and let $\ph\colon F\rtarr F'$ be a natural transformation that commutes with $\SI$. 
Then the full subcategory of $\sA$ whose objects are those $X$ for which $\ph$ is an 
isomorphism is thick, and it is localizing if $F$ and $F'$ preserve coproducts.
\end{prop}
\begin{proof}
Since a retract of an isomorphism is an isomorphism, closure
under retracts is clear. Closure under triangles is immediate from the five 
lemma.  A coproduct of isomorphisms is an isomorphism, so closure under
coproducts holds when $F$ and $F'$ preserve coproducts. 
\end{proof}

\section{The formal isomorphism theorems}

We assume throughout this section that $\sC$ and $\sD$ are closed symmetric monoidal 
categories with compatible triangulations and that $(f^*,f_*)$ is an adjoint 
pair of functors with $f^*$ strong symmetric monoidal and exact.

For the Wirthm\"uller context, we assume in addition that $f^*$ has a left adjoint $f_!$. 
The maps (\ref{proj3})--(\ref{WirVer'}) are then given by (\ref{al}) and Proposition \ref{threeb}.
When 
$${\bar{\pi}}\colon f_!(f^{*}Y\ten X) \rtarr Y\ten f_!X$$
is an isomorphism, the map
$$ \om: f_*X\rtarr f_{!}(X\ten C)$$
is defined. Observe that ${\bar{\pi}}$ is a map between exact left adjoints and that 
${\bar{\pi}}$ and $\om$ commute with $\SI$. The results of the previous section give the 
following conclusion.

\begin{thm}[Formal Wirthm\"uller isomorphism]\label{FWD} 
Let $\sC$ be compactly generated with a generating set $\sG$ such
that ${\bar{\pi}}$ and $\om$ are isomorphisms for $X\in \sG$. Then 
${\bar{\pi}}$ is an isomorphism for all $X\in \sC$. If the objects of 
$\sG$ are dualizable, then $\om$ is an isomorphism for all dualizable $X$. 
If $f^*X$ is compact for $X\in\sG$, then $\om$ is an isomorphism 
for all $X\in \sC$.
\end{thm}

The force of the theorem is that no construction of an inverse to $\om$ 
is required: we need only check that $\om$ is an isomorphism one generating 
object at a time. Proposition \ref{onebyone} explains what is needed for that 
verification. 

For the Grothendieck context, we can use the following basic results of Neeman 
\cite[3.1, 4.1]{Nee} to construct the required right adjoint $f^{!}$ to $f_*$ 
in favorable cases. A main point of Neeman's later monograph \cite{Nee3} and of 
Franke's paper \cite{Fr} is to replace compact generation by a weaker notion that 
makes use of cardinality considerations familiar from the theory of Bousfield localization
in algebraic topology. 

\begin{thm}[Triangulated Brown representability theorem]\label{BR} Let $\sA$ be 
a compactly generated triangulated category. A functor $H\colon \sA^{op}\rtarr \sA b$ that
takes distinguished triangles to long exact sequences and converts coproducts
to products is representable.
\end{thm}

\begin{thm}[Triangulated adjoint functor theorem]\label{TAFT} 
Let $\sA$ be a compactly generated triangulated category and $\sB$ be any
triangulated category. An exact functor $F\colon \sA \rtarr \sB$ that
preserves coproducts has a right adjoint $G$.
\end{thm}
\begin{proof} Take $G(Y)$ to be the object that represents the functor $\sB(F(-), Y)$.
\end{proof} 

The map
$${\pi}\colon Y\ten f_*X \rtarr f_*(f^*Y\ten X)$$
of (\ref{recipy1}) commutes with $\SI$. When ${\pi}$ is an isomorphism, 
$$\ph\colon f^*Y\ten f^{!}Z \rtarr f^{!}(Y\ten Z)$$
is defined and commutes with $\SI$. We 
obtain the following conclusion.

\begin{thm}[Formal Grothendieck isomorphism]\label{FGD} 
Let $\sD$ be compactly generated 
with a generating set $\sG$ such that $f^*Y$ is compact and
${\pi}$ is an isomorphism for $Y\in\sG$. Then $f_*$ has a 
right adjoint $f^{!}$, ${\pi}$ is an isomorphism for all $Y\in\sD$, 
and $\ph$ is an isomorphism for all dualizable $Y$. If the functor $f^{!}$ 
preserves coproducts, then $\ph$ is an isomorphism for all $Y\in\sD$.
\end{thm}
\begin{proof}
As a right adjoint of an exact functor, $f_*$ is exact by Lemma \ref{adexact},
and it preserves coproducts by Lemma \ref{preserve}. Thus $f^{!}$ exists by 
Theorem \ref{TAFT}. Now ${\pi}$ is an isomorphism for all $Y$ by
Proposition \ref{wrap}, $\ph$ is an isomorphism for dualizable $Y$ 
by Proposition \ref{dualGr}, and the last statement holds by Propositions \ref{thickgen}
and \ref{wrap}.
\end{proof}

When $f^{!}$ is obtained abstractly from Brown representability, the
only sensible way to check that it preserves coproducts is to appeal to
Lemma \ref{preserve}, requiring $\sC$ to be compactly generated and $f_*X$ 
to be compact when $X$ is in the generating set. 

For the Verdier-Grothendieck context, we assume that we have a second 
adjunction $(f_!,f^!)$, with $f_!$ exact. We also assume given a map 
$$\hat{\pi}\colon Y\ten f_!X \iso f_!(f^*Y\ten X)$$
that commutes with $\SI$. When $\hat{\pi}$ is an isomorphism, the map
$$\ph\colon f^*Y\ten f^{!}Z \rtarr f^{!}(Y\ten Z)$$
is defined and commutes with $\SI$. Using Proposition \ref{dualGr} and the 
results of the previous section, we obtain the following conclusion.

\begin{thm}[Formal Verdier isomorphism]\label{FVD} 
Let $\sD$ be compactly generated with a generating set $\sG$ such that 
$f^*Y$ is compact and $\hat{\pi}$ is an isomorphism for $Y\in\sG$.
Then $\hat{\pi}$ is an isomorphism for all $Y\in\sD$, and
$\ph$ is an isomorphism for all dualizable $Y$. If the functor $f^{!}$ 
preserves coproducts, then $\ph$ is an isomorphism for all $Y\in\sD$.
\end{thm}

\begin{rem} In many cases, one can construct a more explicit right 
adjoint $f^{!}_{0}$ from some subcategory $\sD_0$ of $\sD$ to some 
subcategory $\sC_0$ of $\sC$, as in Remark \ref{dcontext}. In such 
cases we can combine approaches. Indeed, assume that we have an adjoint 
pair $(f_{!},f^{!}_{0})$ on full subcategories $\sC_0$ and $\sD_0$ 
such that objects isomorphic to objects in $\sC_0$ (or $\sD_0$) are in 
$\sC_0$ (or $\sD_0$). Then, by the uniqueness of adjoints, the right 
adjoint $f^{!}$ to $f_{!}$ given by Brown representability restricts on 
$\sD_0$ to a functor with values in $\sC_0$ that is isomorphic to the 
explicitly constructed functor $f^{!}_{0}$. That is, the right adjoint 
given by Brown representability can be viewed as an extension of the
functor $f^{!}_{0}$ to all of $\sD$. This allows quotation of 
Proposition \ref{three} or \ref{threea} for the construction and 
comparison of the natural maps (\ref{proj2})--(\ref{Ver'}).
\end{rem}

We give an elementary example and then some remarks on the proofs of the results that we 
have quoted from the literature, none of which are difficult.

\begin{exmp}\label{trivex2}
Return to Example \ref{trivex}, but assume further that $\sC$ is a compactly 
generated triangulated category. Here the formal Verdier duality theorem says 
that $\ph=\nu: Y\ten\Hom(C,Z)\rtarr \Hom(C,Y\ten Z)$ is an isomorphism if and 
only if the functor $\Hom (C,-)$ preserves coproducts. That is, an object $C$ 
is dualizable if and only if $\Hom (C,-)$ preserves coproducts.
\end{exmp}

\begin{rem}\label{rkpfs} 
Clearly Theorem \ref{TAFT} is a direct consequence of Theorem \ref{BR}.
In turn, Theorem \ref{BR} is essentially a special case of Brown's original categorical 
representation theorem \cite{Brown}. Neeman's self-contained proof closely 
parallels Brown's argument. The first statement of Proposition \ref{thickgen} 
is used as a lemma in the proof, but it is also a special case. To see
this, let $\sB$ be the localizing subcategory of $\sA$ generated by $\sG$. Then, applied 
to the functor $\sA(-,X)$ on $\sB$ for $X\in\sA$, the representability theorem gives an 
object $Y\in\sB$ and an isomorphism $Y\iso X$ in $\sA$. The second part of 
Proposition \ref{thickgen} is intuitively clear, since objects in $\sA$ not in the thick 
subcategory generated by $\sG$ must involve infinite coproducts, and these will be neither 
dualizable nor compact. The formal proof in \cite{HPS} starts from Example \ref{trivex2}, 
which effectively ties together dualizability and compactness.
\end{rem}

\end{document}